\documentclass[11pt,a4paper]{amsart}
\usepackage[capitalise,noabbrev,nameinlink]{cleveref}

\usepackage{amssymb,verbatim}
\usepackage{amsfonts}
\usepackage[utf8]{inputenc}
\usepackage[mathscr]{euscript}
\usepackage{enumerate}
\usepackage[normalem]{ulem}
\usepackage[
margin=1.5in,
includefoot,
footskip=30pt,
]{geometry} 
\usepackage{mathtools}
\usepackage[dvipsnames]{xcolor}
\usepackage[utf8]{inputenc}
\usepackage[colorinlistoftodos,prependcaption,textsize=tiny]{todonotes}

\newtheorem{thm}[equation]{Theorem}
\newtheorem{lem}[equation]{Lemma}
\newtheorem{prop}[equation]{Proposition}
\newtheorem{cor}[equation]{Corollary}

\theoremstyle{definition}
\newtheorem{dfn}[equation]{Definition}
\newtheorem{definition}[equation]{Definition}

\theoremstyle{remark}
\newtheorem{remark}[equation]{Remark}

\numberwithin{equation}{section}

\newcommand{\N}{\mathbb{N}}
\newcommand{\Z}{\mathbb{Z}}

\newcommand {\G}{\mathcal{G}}
\newcommand {\F}{\mathcal{F}}
\newcommand {\GG}{\mathcal{G}}

\newcommand {\al}{\alpha}

\newcommand{\ann}{\mathrm{ann}}

\title{Characterizing Rickart and Baer Ultragraph Leavitt Path Algebras}
\author{Mitchell Jubeir}
\address[Mitchell Jubeir]{Department of Mathematics, Dartmouth College, Hanover, NH 03755}
\address{Current address: Department of Mathematics, University of California, Santa Barbara, 552 University Rd, Isla Vista, CA 93117}
\email{mitchell\_jubeir@ucsb.edu}
\author{Daniel W. van Wyk}
\address[Daniel W. van Wyk]{Department of Mathematics, Fairfield University, Fairfield, CT 06825 USA, and  Department of Mathematics and Applied Mathematics, University of the Free State, Park West, Bloemfontein, 9301, South Africa.}
\email{dvanwyk@fairfield.edu}
\date{\today}

\begin{document}

\maketitle

\begin{abstract}
    We characterize ultragraph Leavitt path algebras that are Rickart, locally Rickart, graded Rickart, and graded Rickart *-rings. We also characterize ultragraph Leavitt path algebras that are Baer, locally Baer, graded Baer, Baer *-rings, and combinations of these. These characterizations build on and generalize the work of Hazrat and Va\v s on Leavitt path algebras over fields to ultragraph Leavitt path algebras over semi-simple commutative unital rings. 
\end{abstract}
\section{Introduction}

Given a directed graph $E$ and a unital commutative ring $R$, there is a standard process by which to associate an $R$-algebra $L_R(E)$ with $E$, called the Leavitt path algebra of $E$; for example, see  \cite{LeaviittBook} and \cite{Tom11}. These algebras have garnered much attention over the past two decades. For $R=\mathbb{C}$, we may view $L_\mathbb{C}(E)$ as a purely algebraic analog of the  graph $C^*$-algebra $C^*(E)$ of $E$. 

A ring $R$ is a left Baer ring if the left annihilator of any subset of $R$ is a left ideal generated by an idempotent. It is a left Rickart ring if the same condition is satisfied for all singleton sets in $R$. Variations of these notions have been studied when the ring's structure allows for it; for example, if $R$ has an involution or $R$ is a graded ring. 

In \cite{HaVa18}, Hazrat and Va\v s characterize Leavitt path algebras over fields that are Rickart and Baer, including variations such as being locally Rickart, graded Rickart, and graded Rickart *, and similar variations of Baer rings. At least part of the appeal to characterize such properties for $L_R(E)$ lies in that the combinatorial structure of $E$ provides concrete and tangible ways to express abstract algebraic properties of $L_R(E)$.

Building on the work of Hazrat and Va\v s, our goal in the is paper is to extend their results in two directions: (i) we generalize their characterizations to ultragraph Leavitt path algebras, and (ii) generalize the base field to a unital commutative semisimple ring. In particular, our characterizations also apply to algebraic Exel-Laca algebras.

Ultragraphs are similar to graphs but differ in that each edge in an ultragraph may terminate in multiple vertices, as opposed to a single range vertex in the case of graphs. Tomforde introduced ultragraphs in  \cite{Tom03} as a common framework for studying graph $C^*$-algebras and Exel-Laca algebras. It is well-known that the classes of graph $C^*$-algebras and Exel-Laca algebras do not coincide, and neither is one a subclass of the other. However, every ultragraph $C^*$-algebra is (strongly) Morita equivalent to a graph $C^*$-algebra (\cite{KMST10}). 

In \cite{ImaPouLar20}, the authors introduce ultragraph Leavitt path algebras. Every ultragraph Leavitt path algebra is Morita equivalent to a graph Leavitt path algebra (\cite{Fir20}). Therefore, these algebras share any Morita invariant properties such as module or representation theoretic properties. Hence, the study of such properties in ultragraph Leavitt algebras may be confined to the simpler structure of graphs and their Leavitt path algebras. 

Beyond Morita equivalence, however, there are ultragraph Leavitt path algebras that are not isomorphic to any graph algebra; see \cite[Example 5.11]{ImaPouLar20} or \cite[Proposition 2.7]{GonNam23}. It is, therefore, meaningful to study non-Morita invariant properties of ultragraph Leavitt path algebras as generalizations of graph Leavitt path algebras and (algebraic) Exel-Laca algebras (\cite[Theorem 5.7]{ImaPouLar20}).

The properties of being a Baer or Rickart ring are examples of non-Morita invariant properties. A simple example that illustrates this is the ultragraph consisting of an infinite set of vertices $\{v_i: i\in \mathbb{Z}\}$, and an infinite set of edges $\{e_i:i\in \mathbb{Z}\}$ such that $s(e_i)=v_i$ and $r(e_i)=\{ v_{j}: j\geq i\}$. In this case, the ultragraph Leavitt path algebra $L_\mathbb{C}(G)$ is a unital ring, but the graph Leavitt path algebra $L_\mathbb{C}(E)$, which is Morita equivalent to $L_\mathbb{C}(G)$, is a non-unital ring. Therefore, $L_\mathbb{C}(G)$ is a Rickart ring (see Theorem~\ref{thm: Rickart}), but $L_\mathbb{C}(E)$ cannot be Rickart, since all Rickart rings are unital.

We restrict the base ring to a semisimple unital commutative ring to extend Hazrat and Va\v s' characterization from Leavitt path algebras over fields to Leavitt path algebras over rings. If the coefficient ring is not at least Rickart, then the Leavitt path algebra may not be Rickart either (see Remark~\ref{rem:need_semisimple}). We consider semisimple unital commutative rings as the coefficient ring in this paper. Then the coefficient ring decomposes into a finite direct product of fields. Using this fact, we show that a Leavitt path algebra over a semisimple unital commutative ring can be decomposed into a finite direct product of Leavitt path algebras over fields. This decomposition, the Morita equivalence of ultragraph Leavitt path algebras to graph Leavitt path algebras, and Hazrat and Va\u s' characterizations in \cite{HaVa18} are our primary tools for finding various Baer and Rickart-related characterizations for ultragraph Leavitt path algebras.

The structure of the paper is as follows. Section~\ref{s:prelim} contains relevant terminology and definitions of graphs, ultragraphs, and their Leavitt path algebras and an outline of the graph constructed from an ultragraph whose Leavitt path algebra is Morita equivalent to the ultragraph Leavitt path algebra. Section~\ref{s:baer_rickart_rings} includes background on Rickart and Baer rings, including local, graded, and involutive ring versions of these properties. In Section~\ref{s:semi-simple_decompoistion}, we show that an ultragraph Leavitt path algebra over a semisimple ring $R$ can be decomposed into a finite direct product of Leavitt path algebras over fields. This decomposition is essential in our work as it enables us to leverage the Baer and Rickart characterizations of graphs Leavitt path algebras. Section~\ref{s:rickart} and Section~\ref{s:baer} contain our Rickart and Baer characterizations, respectively. 

All rings in this paper are assumed to be associative.

\section{Preliminaries} \label{s:prelim}
In this section we recall the definition of an ultragraph and show how to associate a graph with any ultragraph. We also show that certain combinatorial properties of an ultragraph are inherited by its associated graph.  

\subsection{Graphs and Ultragraphs}
A \emph{(directed) graph} $E=(E^0,E^1,r,s)$ consists of a non-empty countable set  $E^0$ of \emph{vertices}, a countable set $E^1$ of \emph{edges}, and \emph{range} and \emph{source} functions $r,s:E^1\to E^0$. 

A \emph{path of length $n$} in a graph $E$ is a sequence $\lambda=e_1e_2\ldots e_n$ of edges such that $r(e_i)=s(e_{i+1})$ for $i=1,\ldots,n-1$. The \emph{length} of $\lambda$ is $|\lambda|=n$, and vertices are viewed as paths of length $0$. Let $E^n$ denote the set of all paths of length $n$, and let $E^{\ast}=\cup_{n\geq 0}E^n$. Similarly, we define a \emph{path of infinite length} (or an \emph{infinite path}) as an infinite sequence $\lambda=e_1e_2\ldots$ of edges such that $r(e_i)=s(e_{i+1})$ for $i\geq 1$; for such a path, we write $|\lambda|=\infty$ and we let $E^{\infty}$ denote the set of all infinite paths. A vertex $v\in E^0$ is an \textit{infinite emitter} if $|s^{-1}(v)|=\infty$, a \textit{sink} if $|s^{-1}(v)|=0$, and a \textit{source} if $|r^{-1}(v)|=0$. A vertex that is neither a sink nor an infinite emitter is a \textit{regular vertex}.  Let $E_{reg}^0$ denote the set of regular vertices. The graph $E$ is \textit{row-finite} if $E$ has no infinite emitters.

Following the terminology of Hazrat and Va\v{s} in \cite{HaVa18}, a path $e_1\cdots e_n$ in $E$ is \textit{closed} if $s(e_1)=r(e_n)$. If $\mu=e_1\cdots e_n$ is a closed path in $E$ such that $s(e_i)\neq s(e_j)$ for all $i\neq j$, then $\mu$ is a \textit{cycle}. A cycle of length one is a \textit{loop}. A graph $E$ is said to be a \textit{no-exit} graph if $s^{-1}(v)$ has just one element for every vertex $v$ appearing in a cycle.

An infinite path in $E$ is an \textit{infinite sink} if it has no cycles or exits.  An infinite path $e_1e_2\cdots$ in $E$  \textit{ends in a sink or a cycle} if there is a positive integer $n$ such that the path $e_ne_{n+1\cdots}$ is either an infinite sink or equal to a cycle with positive length.

\begin{dfn}
    An \emph{ultragraph} $\mathcal{G}=(G^0, \mathcal{G}^1, r,s)$ consists of a countable set $G^0$ of vertices, a countable set $\mathcal{G}^1$ of edges, a map $s:\mathcal{G}^1 \to G^0$, and a map $r:\mathcal{G}^1 \to \mathcal{P}(G^0)\setminus \{\emptyset\}$, where $\mathcal{P}(G^0)$ denotes the power set of $G^0$.
\end{dfn}

\begin{dfn}
	Let $\mathcal{G}$ be an ultragraph. Define $\mathcal{G}^0$ to be the smallest subset of $\mathcal{P}({G^0})$ that contains $\{v\}$ for all $v\in G^0$, contains $r(e)$ for all $e\in \mathcal{G}^1$, contains $\emptyset$, and is closed under finite unions, finite intersections and relative complements. Elements of $\G^0$ are called \emph{generalized vertices}. 
\end{dfn}

Every graph $E$ is an ultragraph, where, for every $e\in \mathcal{G}^1$, $r(e)=\{v\}$ for some $v\in G^1$. 

Let $\mathcal{G}$ be an ultragraph. A \textit{finite path} in $\mathcal{G}$ is either an element of $\mathcal{G}%
^{0}$ or a sequence of edges $\alpha = \alpha_{1}\ldots \alpha_{k}$ where
$s\left(  \alpha_{i+1}\right)  \in r\left(  \alpha_{i}\right)  $ for $1\leq i\leq k-1$. The \emph{length} of $\alpha$ is $\left|  \alpha\right| =k $. Elements  $A\in\mathcal{G}^{0}$ are considered as paths of length zero. The set of
finite paths in $\mathcal{G}$ is denoted by $\mathcal{G}^{\ast}$. Extend the range and source maps to $\mathcal{G}^{\ast}$ by declaring that $r\left(  \alpha\right)  =r\left(  \alpha_{k}\right) $ and
$s\left(  \alpha\right)  =s\left(  \alpha_{1}\right)$. For $A\in\mathcal{G}^{0}$,
we set $r\left(  A\right)  =A=s\left(  A\right)  $.  An \textit{infinite path} in $\mathcal{G}$ is an infinite sequence of edges $\alpha=\alpha_{1}\alpha_{2}\ldots$, where
$s\left(  \alpha_{i+1}\right)  \in r\left(  \alpha_{i}\right)  $ for all $i$. The set of
infinite paths  in $\mathcal{G}$ is denoted by $
\G^{\infty}$. The length $\left|  \gamma\right|  $ of $\gamma\in\G^{\infty}$ is defined to be $\infty$. A vertex $v$ in $G^0$ is
called a \emph{sink} if $\left|  s^{-1}\left(  v\right)  \right|  =0$, it is
called an \emph{infinite emitter} if $\left|  s^{-1}\left(  v\right)  \right|
=\infty$, and it is called a \emph{source} if $v\notin r(e)$ for all $e\in\GG^1$. A vertex $v$ is \textit{regular} if it is neither a sink  nor an infinite emitter. Let $G_{reg}^0$ denote the set of regular vertices.

A finite path $\al\in \GG^*$ with $|\al|>0$ is \textit{closed} if $s(\al)\in r(\al)$. 
We say $\al$ is a closed path based at $A\in \GG^0$ if $s(\al)\in A$. 
If $\al=\al_1\ldots\al_n$, then $\al$  is a \textit{cycle} if $s(\al_i)\neq s(\al_j)$ for  $i\neq j$; that is, the closed path does not pass through any $s(\al_i)$ multiple times. A cycle of length one is a \textit{loop}.  An \textit{exit} of a path $\al=\al_1\ldots\al_n$ (including $n=\infty$) is either of the following:
\begin{enumerate}
    \item an edge $e\in \GG^1$ such that there exists an $i$ for which $s(e)\in r(\al_i)$ and $e\neq  \al_{i+1}$,
    \item a sink $w$ such that $w \in r(\al_i)$ for some $i$.
\end{enumerate}
We say $\G$ is a \textit{no-exit} ultragraph if no cycle in $\G$ has an exit; that is, for every cycle $\al=\al_1\ldots\al_n$, we have that 
$|s^{-1}(s(\alpha_i))|=|r(\alpha_i)|=1$ for $1\leq i\leq n$.
An infinite path $\alpha=e_{i_1}e_{i_2}\cdots \in \G$ is an \textit{infinite sink} if it has no cycles or exits, that is, $\alpha$ does not contain any cycles and $|s^{-1}(s(e_{i_j}))|=1$ and $|r(e_{i_j})|=1$ for all $j\in \N$. An infinite path $\alpha_1\alpha_2\cdots$ in $E$  \textit{ends in a sink or a cycle} if there is a positive integer $n$ such that the path $\alpha_n\alpha_{n+1\cdots}$ is either an infinite sink or equal to $\gamma^\infty$ for some cycle $\gamma$ with positive length. We say $\G$ is \textit{row-finite} if no vertex $v\in G^0$ is an infinite emitter and $|r(e)|<\infty$ for all $e\in \G^1$. 

\begin{remark}
    We caution that the terminology introduced above is not standard in the literature. A closed path, as defined above, is a loop in \cite{CGvW21} and \cite{CGvW21:mjm}, and a loop above is a simple loop in \cite{CGvW21:mjm}. See also \cite[Proposition 3.7]{CGvW21:mjm}, where isolated points in ultragraph groupoids are characterized in terms of wandering, eventually periodic and semi-tail paths in the ultragraph. These notions are related to that of `infinite sinks' defined above. Since our work builds on \cite{HaVa18}, we opted to use the same terminology as in \cite{HaVa18}.
\end{remark}

\subsection{Graphs Construced from Ultragraphs} \label{construction}
In \cite[Section 3]{KMST10}, a graph is associated with an ultragraph such that the graph $C^*$-algebra is Morita equivalent to the ultragraph $C^*$-algebra. The same graph construction is used in \cite{Fir20} to prove that every ultragraph Leavitt path algebra is Morita equivalent to a graph Leavitt path algebra. In this section, we briefly review the graph construction and refer the reader to \cite{KMST10} for more details and an illuminating example.

Given an ultragraph $\G$, we construct a graph $E$ from $\G$ in the following manner. Begin by fixing an enumeration $\{e_1,e_2,\ldots\}$ of $\G^1$. Define a word $\omega$ to be an element of $\{0,1\}^n$. Let $|\omega|$ denote the length of the word, and let $\omega_i$ denote the $i$th coordinate (or letter) of $\omega$. Also, let $(\omega,0)$ and $(\omega,1)$ denote the word formed by adding a zero and a one, respectively, at the end of $\omega$. For $m\in\N$ and $m\leq n$, if $\omega\in \{0,1\}^{n}$, then we define $\omega|_{m}=(\omega_1,\ldots,\omega_m)\in \{0,1\}^{m}$. The range of a word $\omega$ is the set $r(\omega)=\bigcap_{i: \omega_i =1} r(e_i) \setminus \bigcup_{i: \omega_i =0} r(e_i)$. 

Let $\Delta_n$ be the set of all words of length $n$ such that $|r(\omega)|=\infty$, and let 
\begin{equation}\label{eq:delta_dfn}
    \Delta=\displaystyle \bigsqcup_{n\in\N } \Delta_n.
\end{equation}
Let $\Gamma_0$ denote all words in $\Delta$ that are a string of zeros followed by a single one, and let $\Gamma_+=\Delta \setminus \Gamma_0$. Let $W_+=\bigcup_{\omega \in \Delta} r(\omega)\subset G^0$, and let $W_\infty$ denote the set of all vertices $v\in W_+$ such that $v\in r(\omega)$ for infinitely many $\omega \in \Delta$. 

By \cite[Lemma 3.7]{KMST10}, there exists a function $\sigma : W_+ \to \Delta$ defined as follows. If $v\notin W_\infty$, then $\sigma(v)$ is the longest word in $\Delta$ such that $v\in r(\omega)$. Fix an enumeration of $W_\infty$. For $v_k\in W_\infty$, let $N_k$ denote the set $\{\omega\in\Delta :v_k\in r(\omega),\,|\omega|>k\}$.  
Let $n$ denote the length of the shortest word in $N_k$. Define $\sigma(v_k)$ to be the unique $\omega\in \Delta_n$ such that $v_k\in r(\omega)$.
Notice that $v\in r(\sigma(v))$ for each $v\in W_+$, and that $\sigma^{-1} (\omega)$ is finite, and possibly empty, for each $\omega \in \Delta$. Extend $\sigma$ to a function $\sigma: G^0\to \Delta\cup \emptyset$ by putting $\sigma(v)=\emptyset$ for all $v\notin W_+$.

Next, define $X(e_n)=\{v\in r(e_n) : |\sigma(v)|<n\} \sqcup \{\omega \in \Delta_n : \omega_n=1\}$. 
For each $n\in \N$, the set $X(e_n)$ is nonempty and finite (\cite[Lemma 3.11]{KMST10}).

We can now define a graph $E=(E^0,E^1,r,s)$ constructed from $\G$. Let 
\begin{align*}
    E^0 & = G^0 \sqcup \Delta, \\
    E^1 & = \{\bar{x}: x\in W_+\sqcup \Gamma_+\} \sqcup \{(e_n,x): e_n\in \G^1, x\in X(e_n)\}, \\
    r_E(\bar{x}) & = x, \, \, r_E((e_n,x)) = x, \text{ and } \\
    s_E(\bar{v}) &= \sigma(v), \,\, s_E(\bar{\omega})=\omega|_{|\omega|-1}, \,\, s_E((e_n,x)) = s(e_n) \\
\end{align*}
We refer to $E$ as  \textit{the graph associated with} $\G$. The construction of $E$ relies on a fixed enumeration of $\G^1$ and a choice for $\sigma$. However, the Morita equivalence between an ultragraph algebra and its associated graph's algebra is independent of these choices (\cite[Section 5]{KMST10}, \cite[Section 9]{Fir20}). 

A key result is a correspondence between paths in $E$ and paths in $\G$.
\begin{thm}\cite[Lemma 4.9]{KMST10}\label{thm:pathcorrespondance}
Let $\G $ be an ultragraph and $E$ the graph associated with $\G$, and fix $v, w\in G^0$. There is a bijection between paths in $E$ from $v$ to $w$, and paths in $\G$ beginning at $v$ whose ranges contain $w$. 
\end{thm} 

We prove some preliminary lemmas about the combinatorial structure of the graph associated with an ultragraph. Fix an ultragraph $\G$ and let $E$ be its associated graph for any enumeration of $\G^1$ and choice of $\sigma$. Let $F^1=\{\overline{x} :x\in W_+\sqcup \Gamma_+ \} \subset E^1$; that is, $F^1$ is the set of all edges in $E^1$ starting at some $\omega\in\Delta$. Let $F^* =E^0\sqcup F^1$.

\begin {lem}\label{lem: row-finite} The following hold.
\begin{enumerate}[(i)]
    \item $E$ is row-finite if and only if $\G$ is without infinite emitters.
    \item $E$ is without sinks if and only if $\G$ is without sinks.
\end{enumerate}
\end {lem}
\begin{proof}
By \cite[Proposition 3.14]{KMST10}, we have that $E^0_{reg} = G^0_{reg}\sqcup \Delta$. Thus, for all vertices $\overline{v}\in \Delta$, $s^{-1}(\overline{v})$ is nonempty and finite. And, for any $v\in G^0$,  $s_E^{-1}(v)$ is finite if and only if $s^{-1}(v)$ is finite, and $s_E^{-1}(v)$ is nonempty if and only if $s^{-1}(v)$ is nonempty. Hence $E$ is row-finite if and only if $\G$ is without infinite emitters and $E$ is without sinks if and only if $\G$ is without sinks.
\end{proof}

\begin{lem}\label{l:delta_empty_row_finite}
  Fix an enumeration $\{e_1,e_2,\ldots\}$ of $\G^1$. Then, $|r(e_i)|<\infty$ for all $i \in \N$ if and only if $\Delta=\emptyset$.  
\end{lem}
\begin{proof}
    Suppose $|r(e_i)|<\infty$ for all $i\in \N$. Then $|\cap_{j=1}^nr(e_{i_j})|<\infty$ for all $n\in \N$ and $e_{i_j}\in \G^1$. Thus,  $|r(\omega)|<\infty$ for all $\omega\in\{0,1\}^n$ and for all $n\in \N$. Hence, $\Delta=\emptyset$.

    Suppose $\Delta=\emptyset$. To the contrary, suppose that there is an $n\in \N$ such that $|r(e_n)|=\infty$. Let $i=\min\{ j\in\N: |r(e_j)|=\infty\}$. Consider the word $\omega=(0^{i-1},1) \in\{0,1\}^i$ consisting of a string of $i-1$ zeros followed by a $1$. Then 
    \[|r(\omega)|=|r(e_i)|=\infty, \]
    which contradicts that $\Delta=\emptyset$. Hence $|r(e_i)|<\infty$ for all $i \in \N$.
\end{proof}

\begin{prop}\label{prop:Str_row_finite_equiv}
    $\G$ is row-finite if and only if $E$ is row-finite and $\Delta=\emptyset$. 
\end{prop}
\begin{proof}
Follows from Lemma~\ref{lem: row-finite} and Lemma~\ref{l:delta_empty_row_finite}
\end{proof}

\begin{lem}\label{lem: inf path sink/cycles} 
$E$ is a row-finite, no-exit graph in which every infinite path ends in a sink or a cycle if and only if $\G$ is row-finite, no-exit, and all infinite paths end in sinks or cycles. 
\end {lem}
\begin{proof}

Suppose  $E$ is a row-finite, no-exit graph in which every infinite path ends in a sink or a cycle. We claim that $\Delta =\emptyset$. Suppose $\Delta$ is non-empty. Then Lemma~\ref{l:delta_empty_row_finite} implies there exists a word $\omega\in \Delta_k$, for some $k\in\N$. Then, either $(\omega,0)$ or $(\omega,1)$ is a member of $\Delta_{k+1}$ (see \cite[Remark 3.3]{KMST10}) and there is an edge connecting $\omega$ to this word. Thus, there is an infinite path along words. Suppose it starts with $\omega_0 \in \Gamma_0$, then there is an infinite path $\alpha=\alpha_1\alpha_2\cdots$ in $F^*$. Suppose that $\alpha$ ends in a sink (since it cannot end in a cycle). Then there is an $n\in\N$ such that there are no exits after the $n$th term.  Let $\omega^\prime=s(\alpha_i)$ be a word in $\Delta$ for some $i>n$, and suppose the $\omega^\prime\in \Delta_k$. Then, $\{ v\in r(\omega^\prime) : |\sigma(v)|>k\}$ must be infinite since $\sigma^{-1}(\omega)$ is finite for all words $\omega\in\Delta$ and there are finitely many words of length less than $k$. Thus, by \cite[Lemma 4.6]{KMST10}, there is a path from $\omega^\prime$ to each $v\in\{ v\in r(\omega^\prime) : |\sigma(v)|>k\}$, so there must be an exit from the path. Thus, this path must have an exit, and, therefore, if $\Delta$ is non-empty, then there is an infinite path not ending in sink or cycle, contradicting our assumptions. Hence $\Delta=\emptyset$. 

On the other hand, if $\G$ is strongly row finite, then $\Delta=\emptyset$ by Proposition~\ref{prop:Str_row_finite_equiv}. Thus, both the necessary and the sufficient hypotheses imply that $\Delta=\emptyset$. 

Since $\Delta=\emptyset$, it follows that $G^0=E^0$, and then the result follows from the path correspondence given by Theorem~\ref{thm:pathcorrespondance}.
\end{proof}

\subsection{Leavitt Path Algebras of Ultragraphs} \,

Following \cite{ImaPouLar20}, the Leavitt path algebra of an ultragraph is defined as follows.
\begin{dfn}
Let $\mathcal{G}$ be an ultragraph and $R$ a unital commutative ring. The \textit{ultragraph Leavitt path algebra} $L_R(\mathcal{G})$ of $\mathcal{G}$ is the universal $R$-algebra generated by $\{s_e,s_e^*:e\in \mathcal{G}^1\}\cup\{p_A:A\in \mathcal{G}^0\}$ subject to the relations
	\begin{enumerate}
		\item $p_\emptyset=0,  p_Ap_B=p_{A\cap B},  p_{A\cup B}=p_A+p_B-p_{A\cap B}$, for all $A,B\in \mathcal{G}^0$;
		\item $p_{s(e)}s_e=s_ep_{r(e)}=s_e$ and $p_{r(e)}s_e^*=s_e^*p_{s(e)}=s_e^*$ for each $e\in \mathcal{G}^1$;
		\item $s_e^*s_f=\delta_{e,f}p_{r(e)}$ for all $e,f\in \mathcal{G}^1$;
		\item $p_v=\sum\limits_{s(e)=v}s_es_e^*$ whenever $0<\vert s^{-1}(v)\vert< \infty$.
	\end{enumerate}
\end{dfn}

\begin{remark}
In \cite{Tom03}, $\G^0$ is not required to be closed under relative complements, whereas in  \cite{ImaPouLar20} it is. By \cite[Proposition 5.2]{CGvW21}, the Leavitt path algebra of $\G$ with $\G^0$ closed under relative compliments is isomorphic to the Leavitt path algebra of $\G$ without this assumption on $\G^0$. Thus, there is no harm in allowing relative complements here.  
\end{remark} 

If $\alpha=\alpha_1\ldots\alpha_n\in \mathcal{G}^*$, with $\alpha_i\in\mathcal{G}^1$, then $s_\alpha:=s_{\alpha_1}\ldots s_{\alpha_n}$, and $s^*_\alpha:=s^*_{\alpha_n}\ldots s^*_{\alpha_n}$.
By \cite[Theorem 2.9]{ImaPouLar20}
\begin{equation*}
    L_R(\mathcal{G}) = \mathrm{span}_R\{s_\alpha p_A s_\beta^*: \alpha,\beta\in \mathcal{G}^*, A\in \mathcal{G}^0\}.
\end{equation*}
By \cite[Lemma 2.12]{ImaPouLar20}, $L_R(\mathcal{G})$ is unital if and only if $G^0\in \mathcal{G}^0$, with the unit given by $p_{G^0}$. 

Every graph is an ultragraph. If $\G$ is a graph, then the definition above and that of the graph Leavitt path algebra (\cite{LeaviittBook}) coincide. Hence, an ultragraph Leavitt path algebra generalizes graph Leavitt path algebras. 

If the ring $R$ has an involution $\#$ defined on it, then the map $*:rs_e\mapsto r^{\#}s_{e}^*$ extends to an involution on $L_R(\G)$. 

\begin{definition}\label{dfn:grading}
Let $S$ be a ring and $H$ an abelian group. Then, $S$ together with a collection of additive subgroups $\{S_n\}_{n\in H}$ of $S$ such that 
    \begin{enumerate}
        \item $S=\bigoplus_{n\in H}S_n$, and 
        \item $S_mS_n \subset S_{m+n}$ for all $m,n\in H$
    \end{enumerate}
is an \textit{$H$-graded ring}. The subgroup $S_n$ is called the \textit{homogeneous component of degree $n$}, and the collection $\{S_n\}_{n\in H}$ is called an \emph{$H$-grading} of $S$. If $S$ is an $H$-graded ring and a *-ring, then $S$ is an \textit{$H$-graded *-ring} if $S_n^*\subseteq S_{-n}$. If $A$ is an $R$-algebra, then $A$ is a \textit{graded algebra} if $A$ is a graded ring with an $H$-grading $\{A_n\}_{n\in H}$, and each $A_n$ is an $R$-submodule of $A$. 
\end{definition}

\begin{definition}
If $S$ is an $H$-graded ring, an ideal $I\subset S$ is an \emph{$H$-graded ideal} if $I=\bigoplus_{n\in H} (I\cap S_n)$. If $\phi:S\to T$ is a ring homomorphism between $H$-graded rings, then $\phi$ is an \textit{$H$-graded homomorphism} if $\phi(S_n)\subset T_n$ for every $n\in H$.
\end{definition}

We say that an $H$-graded ring $S$ is \textit{graded von Neumann regular} if, for every homogeneous element $x$, there is a homogeneous element $y$ such that $x=xyx$.

Every ultragraph Leavitt path algebra $L_R(\G)$ is a $\Z$-graded ring with homogeneous components given by 
\[
L_R(\G)_n = \mathrm{span}_R \left\{ s_\alpha p_A s_{\beta}^*\in L_R(\G) : \alpha,\beta\in \G^* \text{ and } |\alpha|-|\beta| = n \right\}.
\]

Since $R$ is commutative, it always has an involution (also denoted by *). Thus, any ultragraph Leavitt path algebra over a commutative unital ring $R$ is an involutive algebra, with the involution given by 
\[ \left(\sum_{i=1}^n r_i s_{\alpha_i}p_As_{\beta_i}^*\right)^*= \sum_{i=1}^n r_i^*s_{\beta_i}p_A s_{\alpha_i}^*, 
\] 
where $\alpha_i,\beta_i\in \mathcal{G}^*$ and $A\in\G^0$.

\section{Baer and Rickart rings.} \label{s:baer_rickart_rings}
This section recalls the definitions of Rickart and Baer rings and variations relevant to our goal of characterizing ultragraph algebras with these properties. We follow \cite[Section 1]{HaVa18} fairly closely and refer the reader there for more details.

Let $S$ be a ring.  For a subset $B\subset S$, the \textit{left annihilator of $B$} is the set $\ann_l(B)=\{s\in S: sb=0 \text{ for every } b\in B \}$, which is a left ideal of $S$. The right annihilator of $I$ is $\ann_r(B)=\{s\in S: bs=0 \text{ for every } b\in B \}$, which is a right ideal of $S$.

\textit{Rickart and Baer rings.} 
An \textit{idempotent} in a ring $S$ is an element $e\in S$ such that $e^2=e$. The ring $S$ is a \textit{left Rickart ring} if, for all $x\in S$, $\ann_l(x)$ is generated as a left ideal by an idempotent; that is, $\ann_l(x)=eS$ for some idempotent $e\in S$. A right Rickart ring is defined similarly. If $S$ is both a left and right Rickart ring, then we say $S$ is a \textit{Rickart ring} (or Rickart for short). Rickart rings are necessarily unital. We say $S$ is a \textit{left Baer ring} if, for every subset $B\subset S$, $\ann_l(B)$ is generated as a left ideal by an idempotent. 
Right Baer rings are defined similarly, and a ring that is both a left and a right Baer ring is a Baer ring.

\textit{Rickart and Baer *-rings.}
If $S$ is a *-ring, then $\ann_l(x^*) =(\ann_r(x))^*$. Thus, $S$ is a left Rickart ring if and only if $S$ is a right Rickart ring, and in either case, $S$ is necessarily unital. A projection is an idempotent $p\in S$ such that $p^*=p$. For *-rings, projections play the role that idempotents play in the definitions of Rickart and Baer rings. Hence, we say $S$ is a \textit{Rickart *-ring} if $\ann_l(x)$ is generated as a left ideal by a projection, or equivalently, if $\ann_r(x)$ is generated as a right ideal by a projection. Rickart *-rings are also necessarily unital rings. Similarly, $S$ is a Baer *-ring if $\ann_l(B)$ (or equivalently $\ann_r(B)$) is generated as a left (or right) ideal by a projection.

\textit{Graded Rickart and Baer rings.}
Let $H$ be a group and let $S$ be an $H$-graded ring. Then $S$ is a \textit{graded left Rickart ring} if for any homogeneous element $x\in S$ the ideal $\ann_l(x)$ is generated by a homogeneous idempotent. Graded right Rickart rings are defined similarly. A \textit{graded Rickart ring} is both a graded left and a graded right Rickart ring. If $B\subset S$ consists of homogeneous elements, $\ann_l(B)$ is a graded ideal. If, for every subset $B\subset S$ consisting of homogeneous elements, $\ann_l(B)$ (or equivalently $\ann_r(B)$) is generated by a homogeneous idempotent, then $S$ is a \textit{graded Baer ring}. 

If $S$ is a graded *-ring, then $S$ a \textit{graded Rickart *-ring} if for every homogeneous $x\in S$, the ideal $\ann_l(x)$ (or equivalently $\ann_r(x)$) is generated by a homogeneous projection. If, for all $B\subset S$ consisting of homogeneous elements, $\ann_l(B)$ (or equivalently $\ann_r(B)$) is generated by a homogeneous projection, then $S$ is a \textit{graded Baer *-ring}.

\textit{Locally Rickart and Baer rings.}
Let $S$ be a ring. A set of idempotents $F\subset S$ is a \textit{set of local units if} if for every finite set $\{s_1, \ldots, s_n\}\subset S$, there is an $f\in F$ such that $fs_i =s_if =s_i$ for every $1\leq i\leq n$. If $S$ has a set of local units, the $S$ is \textit{locally unital}. We say $S$ has \textit{enough idempotents} if there exists a set
of nonzero orthogonal idempotents $E$ for which the set $F$ of finite sums of
distinct elements of $E$ is a set of local units for $S$.

The notions of being Rickart or Baer can be extended to non-unital rings as follows. A ring $S$ is \textit{locally  left Rickart} if the corner $eSe$ is left Rickart for any idempotent $e\in S$. Locally right Rickart, locally Rickart, locally (left/right) Baer, locally Rickart *, locally Baer *, and the graded versions are defined analogously.

We end this section with a few lemmas that we use to prove our main results. The following lemma is taken from the discussion following Definition 6 in \cite{HaVa18}. 

\begin{lem}\cite[Section 1.6]{HaVa18} \label{l:local_units_rickart}
    Let $S$ be a locally unital ring with a local set of units $U$. Then $uSu$ is right Rickart, left Rickart or Baer, respectively, for all $u\in U$ if and only if $S$ is locally right Rickart, locally left Rickart or locally Baer, respectively.
    Analogous statements for graded rings, *-rings, and for graded *-rings hold.
\end{lem}

\begin{lem}\label{l:local_units_direct_prod}
    Let $S=\Pi_{i=1}^n S_i$ be a finite direct product of rings $S_i$. Then $S$ is locally unital if and only if each $S_i$ is locally unital.  
\end{lem}
    \begin{proof}
        Suppose $S$ has a set of local units $U$. Let $\pi_i$ denote the projection of $\Pi_{i=1}^n S_i$ onto $S_i$. Then, a simple computation shows that $U_i=\pi_i(U)$ is a set of local units for $S_i$. Hence, every $S_i$ is locally unital.

        Conversely, if every $S_i$ has a set $U_i$ of local units, then $U=\Pi_{i=1}^n U_i$ is a set of local units for $S$.  
    \end{proof}

\begin{lem}\label{l:Rick_Baer_direct_prod}
    Let $S=\Pi_{i=1}^n S_i$ be a finite direct product of *-rings $S_i$. Then 
    \begin{enumerate}[(i)]
        \item $S$ is a Baer ring if and only if each $S_i$ is a Baer ring.
        \item $S$ is a Rickart ring if and only if each $S_i$ is a Rickart ring.
        \item $S$ is a Baer (Rickart, respectively) *-ring if and only if each $S_i$ is a Baer (Rickart, respectively) *-ring.
        \item $S$ is locally Baer (Rickart, respectively) if and only if each $S_i$ is locally Baer (Rickart, respectively).
        \item Suppose in addition that $S$ is a locally unital graded *-ring such that $s=(s_1,\ldots,s_n)\in S$ is homogeneous if and only if each $s_i\in S_i$ is homogeneous. Then $S$ is a graded locally Baer (Rickart, respectively) $*$-ring if and only if each $S_i$ is a graded locally Baer (Rickart, respectively) $*$-ring. 
    \end{enumerate}
\end{lem}
\begin{proof}
    Items (i) and (ii) follow from the observations that if $B =\Pi_{i=1}^n B_i \subset S$, then  
    \begin{equation}\label{eq:ann_direct_prod}
        \ann_l(B) = \Pi_{i}^n \ann_l(B_i),
    \end{equation}
    and that $e=(e_1,\ldots,e_n)\in S$ is an idempotent if and only if $e_i\in S_i$ is an idempotent for every $1\leq i\leq n$.

     (iii) This follows from Equation~(\ref{eq:ann_direct_prod}) and the observation that $e=(e_1,\ldots,e_n)\in S$ is a projection if and only if $e_i\in S_i$ is a projection for every $1\leq i\leq n$.

     (iv) Suppose $S=\Pi_{i=1}^n S_i$ is a locally Baer ring. By Lemma~\ref{l:local_units_direct_prod}, $u=(u_i,\ldots,u_n)\in S$ is a local unit if an only if $u_i\in S_i$ is a local unit, for every $1\leq i\leq n$. Then, by (i) (and (ii), respectively), $uSu = \Pi_{i=1}^n u_iS_iu_i$ is Baer (Rickart, respectively) if and only if each $u_iS_iu_i$ is Baer (Rickart, respectively). It now follows from Lemma~\ref{l:local_units_rickart} that $S$ is locally Baer (locally Rickart, respectively) if and only if each $S_i$ is locally Baer (locally Rickart, respectively).

     (v) If $S$ is a graded *-ring which is Baer (Rickart, respectively) *-ring, then it is a graded Baer (Rickart, respectively) *-ring (see discussion after Definition 4 in \cite{HaVa18}). Hence, (v) follows from (iii). 
     \end{proof}

\section{Leavitt path algebras over semisimple rings} \label{s:semi-simple_decompoistion}
This section shows that an ultragraph Leavitt path algebra over a semisimple unital commutative ring $R$ can be decomposed into a finite direct product of Leavitt path algebras over fields. This enables us to deduce certain properties of Leavitt path algebras over rings from Leavitt path algebras over fields.

\begin{remark} \label{rem:need_semisimple}
We need to assume that the coefficient ring $R$ is at least a Rickart ring, in addition to being unital and commutative. If $R$ is not a Rickart ring, then there is a nonzero $r\in R$ such that $\ann_l(\{r\})$ is a proper ideal that is not generated by an idempotent. Thus, if $E$ is a graph consisting of a single vertex $v$, then $\ann_l(\{rv\}) = \ann_l(\{r\})v \subsetneq L_R(E)$, which cannot be generated by an idempotent.
We will, however, impose a slightly stronger condition on $R$ than being Rickart. We assume $R$ is semisimple, unital and commutative. This implies that $R$ is Rickart, and it allows us to use the direct product decomposition described below.
\end{remark}  

If $R$ is a unital commutative semisimple ring, then it follows from the Wedderburn–Artin Theorem and the commutativity of $R$ that  $R$ is isomorphic to a finite direct product $\Pi_i^n F_i$, where each $F_i$ is a field. Throughout this paper we write $R=\Pi_i^n F_i$ and $r=(r_1, \ldots, r_n)\in R$ without reference to the isomorphism between $R$ and $\Pi_i^n F_i$. Since each $F_i$ may be viewed as subring of $R$, any involution on $R$ restricts to an involution on $F_i$, such that $(r_1^*,\ldots,r_n^*)=(r_1,\ldots,r_n)^*$. 

Conversely, if each $F_i$ is a *-ring, then $(r_1,\ldots,r_n)^*=(r_1^*,\ldots,r_n^*)$ defines an involution on $R$. That is, $R$ and $\Pi_i^n F_i$ are *-isomorphic rings.

Fix a unital commutative semisimple ring $R =\Pi_i^n F_i$, where each $F_i$ is a field. 

\begin{lem} \label{l:direct_prod_grading}
Let $\G$ be an  ultragraph. Then the direct product $\Pi_{i}^n L_{F_i}(\G)$ of ultragraph Leavitt path algebras is an $R$-algebra which is $\mathbb{Z}$-graded with homogeneous components 
\begin{equation*}
    (\Pi_{i}^n L_{F_i}(\G))_k =  \{ (x_1,\ldots,x_n) : x_i\in L_{F_i}(\G)_n \}
\end{equation*}
for $k\in \mathbb{Z}$. Moreover, if $R$ is a *-ring, then $\Pi_{i}^n L_{F_i}(\G)$ is an involutive $R$-algebra, with the involution defined  componentwise $(x_1,\ldots,x_n)^* = (x_1^*, \ldots,x_n^*)$.
\end{lem}
\begin{proof}
    It is not hard to show that $\Pi_{i}^n L_{F_i}(\G)$ is an (involutive) $R$-algebra, with scalar multiplication given by $(r_1,\ldots, r_n)(x_1,\ldots, x_n) = (r_1x_1,\ldots, r_nx_n)$, where $(r_1,\ldots,r_n)\in R$ and $(x_1,\ldots, x_n)\in \Pi_{i}^n L_{F_i}(\G)$; we omit the details.

    We show that $\Pi_{i}^n L_{F_i}(\G)$ is $\Z$-graded. Let $x=(x_1,\ldots,x_n)\in \Pi_{i}^n L_{F_i}(\G)$. Since each $L_{F_i}(\G)$ is $\Z$-graded by \cite[Theorem 2.9]{ImaPouLar20}, we can express each $x_i$ as a sum of homogeneous components (by Definition~\ref{dfn:grading}(1)). Hence, $x_i = \sum_{j\in\Z} r_{i_j}x_{i_j}$, where $x_{i_j}\in L_{F_i}(\G)_j$. Then 
    \begin{eqnarray}
        x &=& \left(\sum_{j\in\Z} r_{1_j}x_{1_j}, \ldots, \sum_{j\in\Z} r_{n_j}x_{n_j}   \right) \nonumber \\
        &=& \sum_{j\in\Z} ( r_{1_j}x_{1_j}, 0,\ldots,0)+ \cdots + \sum_{j\in\Z}(0,\ldots,0, r_{n_j}x_{n_j} ) \nonumber \\
        &=&  \sum_{j\in\Z} (r_{1_j}x_{1_j}, \ldots,r_{n_j}x_{n_j}).\label{eq:1}
    \end{eqnarray}
     Hence $\Pi_{i}^n L_{F_i}(\G) \subseteq \bigoplus_{j\in \Z} \left(\Pi_{i}^n L_{F_i}(\G)_j\right)$. The reverse inclusion is clear, and thus $\Pi_{i}^n L_{F_i}(\G) = \bigoplus_{j\in \Z} \left(\Pi_{i}^n L_{F_i}(\G)_j\right)$. Since, each $L_{F_i}(\G)_m$ is an abelian subgroup of $L_{F_i}(\G)$, it follows that $\Pi_{i}^n L_{F_i}(\G)_j$ is an abelian subgroup of $\Pi_{i}^n L_{F_i}(\G)$. Hence, condition (1) of Definition~\ref{dfn:grading} is satisfied. In addition, $\Pi_{i}^n L_{F_i}(\G) = \bigoplus_{j\in \Z} \left(\Pi_{i}^n L_{F_i}(\G)_j\right)$ implies that 
    \[(\Pi_{i}^n L_{F_i}(\G))_k = \Pi_{i}^n L_{F_i}(\G)_k = \{ (x_1,\ldots,x_n) : x_i\in L_{F_i}(\G)_k \}.\]

    Next suppose that $x=(x_1,\ldots,x_n)\in (\Pi_{i}^n L_{F_i}(\G))_p$ and $y=(y_1,\ldots,y_n)\in (\Pi_{i}^n L_{F_i}(\G))_q$, with $p,q\in\Z$.  Then $xy= (x_1y_1, \ldots, x_ny_n)$. Since each $L_{F_i}(\G)$ is $\Z$-graded, it follows that $x_iy_i \in L_{F_i}(\G)_{p+q}$ for $1\leq i\leq n$. Thus, $xy\in L_R(\G)_{p+q}$ and $\Pi_{i}^n L_{F_i}(\G)$ is a $\Z$-graded ring. 
\end{proof}

\begin{prop}\label{prop:direct_prod_iso}
Let $\G$ be an ultragraph. Then $L_R(\G)$ is $\mathbb{Z}$-graded *-isomorphic to $\Pi_i^n L_{F_{i}}(\G)$. 
\end{prop}
\begin{proof}
    Let $x\in L_R(\G)$. By \cite[Theorem 9.9]{ImaPouLar20}, we may take $x = \sum_{i=1}^m r_is_{\alpha_i}p_{A_i}s_{\beta_i}^*$, where $\alpha_i,\beta_i\in \G^*$, $A_i\in \G^0$, and $r_i\in R$ for $1\leq i\leq m$. For each $1\leq i\leq m$, we have that $r_i = (r_{i,1},\ldots,r_{i,n})$, where $r_{i,j}\in F_j$. Define 
    $\phi: L_R(\G)\rightarrow \Pi_{i=1}^m L_{F_i}(\G)$ by 
    \[ \phi(x) = \left(\sum_{i=1}^m r_{i,1}s_{\alpha_i}p_{A_i}s_{\beta_i}^*,\ldots,\sum_{i=1}^m r_{i,n}s_{\alpha_i}p_{A_i}s_{\beta_i}^*\right).\]
    It is not hard to see that $\phi$ is a well-defined $\mathbb{Z}$-graded homomorphism of $L_R(\G)$ into $\Pi_i^n L_{F_{i}}(\G)$. 

    To see that $\phi$ is surjective, fix an arbitrary \\ $(\sum_{i=1}^{k_1}r_{i,1}s_{\alpha_{i,1}}p_{A_{i,1}}s_{\beta_{i,1}}^*,\ldots,\sum_{i=1}^{k_n}r_{i,n}s_{\alpha_{i,n}}p_{A_{i,n}}s_{\beta_{i,n}}^*)\in \Pi_i^n L_{F_{i}}(\G)$, where $s_{\alpha_{i,j}}, s_{\beta_{i,j}}\in \G^*$ and $p_{A_{i,j}}\in \G^0$. For $j=1,2,\ldots, n$, let $c_{i,j}\in \Pi_{i}^n F_i$ have $r_{i,j}$ in the $j$th coordinate and zeros elsewhere. Then, $\sum_{j=1}^n\sum_{i=1}^{k_j}c_{i,j}s_{\alpha_{i,j}}p_{A_{i,j}}s_{\beta_{i,j}}^*\in L_R(\G)$, and 
    \begin{equation*}
    \phi\left(\sum_{j=1}^n\sum_{i=1}^{k_j}c_{i,j}s_{\alpha_{i,j}}p_{A_{i,j}}s_{\beta_{i,j}}^*\right) = \left(\sum_{i=1}^{k_1}r_{i,1}s_{\alpha_{i,1}}p_{A_{i,1}}s_{\beta_{i,1}}^*,\ldots,\sum_{i=1}^{k_n}r_{i,n}s_{\alpha_{i,n}}p_{A_{i,n}}s_{\beta_{i,n}}^*\right). 
    \end{equation*}
    Hence, $\phi$ is surjective. 

    To see that $\phi$ is injective, let $r=(r_1,\dots,r_n)\in R\setminus \{0\}$ and $A\in \G^0$ nonempty. Then, since $r_i\neq 0$ for at least one $i\in\{1,\ldots,n\}$, it follows from \cite[Theorem 2.10]{ImaPouLar20} that $rp_{A} \neq 0$. Thus, $\phi(rp_A) \neq 0$ for all $r\in R\setminus\{0\}$ and for all nonempty $A\in\G^0$. Hence, $\phi$ is injective by the Graded Uniqueness Theorem (\cite[Theorem 2.14]{ImaPouLar20}), and thus a graded isomorphism. 

    It is straightforward to show that this isomorphism preserves the involution and is, therefore, a *-isomorphism. 
\end{proof}

\begin{cor}\label{cor:graded_regular}
    The $R$-algebra $L_R(\G)$ is graded von Neumann regular if and only if and $L_{F_{i}}(\G)$ is graded von Neumann regular for every $1\leq i\leq n$.
\end{cor}
\begin{proof}
    The result follows from the fact that $x=(x_1,\ldots,x_n)\in L_R(\G)$ is homogeneous if and only if each $x_i\in L_{F_i}(\G)$ is homogeneous.
\end{proof}

It is known that if $\{A_i\}_{i=1}^n$ is a collection of $R$-algebras such that each $A_i$ is hereditary, then $\Pi_{i=1}^n A_i$ is hereditary (see, for example, Proposition a of Section 6.1 in \cite{Pierce:AssocAlgs} and \cite[Examples 2.32(g)]{Lam_1}). The key is that each component algebra has the same coefficient ring. We use this fact to show that every ultragraph Leavitt path algebra is a hereditary (and thus semihereditary) ring.

\begin{prop}\label{prop:hereditary}
    Let $\G$ be an ultragraph and $R$ a commutative semi-simple ring with unity. Then $L_R(\G)$ is a (semi)hereditary ring.
\end{prop}
\begin{proof}
    First, assume that $\G$ is an ultragraph without singular vertices. Then, by \cite[Theorem 9.6]{Fir20}, $L_R(\G)$ is Morita equivalent to $L_R(E)$, where $E$ is the graph associated with $\G$. Since $R$ is a semisimple commutative ring, there are finitely many fields $F_1, \ldots, F_n$ such  that $L_R(E) = \Pi_{i=1}^{n}L_{F_i}(E)$. By \cite[Theorem 3.7]{AG12}, each $L_{F_i}(E)$ is a hereditary ring. Each $L_{F_i}(E)$ may be viewed as a subalgebra of $L_R(E) = \Pi_{i=1}^{n}L_{F_i}(E)$, and thus as an $R$-algebra, where, for $(r_1, \ldots, r_n)\in R$ and $x\in L_{F_i}(E)$, we have that 
    \[(r_1, \ldots, r_n)x = r_ix.\] 
    Thus, \cite[Section 6.1, Proposition a]{Pierce:AssocAlgs} implies that $L_R(E) = \Pi_{i=1}^{n}L_{F_i}(E)$ is a hereditary ring. 
    
    Now suppose that $\G$ is a general ultragraph, possibly with singular vertices. Then $L_R(\G)$ is Morita equivalent to $L_R(\tilde{\G})$, where $\tilde{\G}$ is the desingulirazation of $\G$, \cite[Theorem 10.6]{Fir20}. Since, $\tilde{\G}$ has no singular vertices, it is hereditary by the first part of the proof. Since hereditary is a Morita invariant property, it follows that $L_R(\G)$ is a hereditary ring.
    
\end{proof}

\section{Rickart ultragraph Leavitt path algebras}\label{s:rickart}
In this section, we characterize ultragraph Leavitt path algebras which are (locally) Rickart, graded Rickart, and graded Rickart *-rings, generalizing analogous characterizations for Leavitt path algebras of graphs (\cite[Proposition 13]{HaVa18}). Throughout this section, $R$ denotes a unital commutative semisimple *-ring.

\begin{definition} 
Suppose $S$ is a *-ring and $s_1,\ldots,s_n\in S$ for an arbitrary $n\in\N$. The involution on $S$ is  \emph{positive definite} if $\sum_{i=1}^n s_i s_i^* =0$ implies that $s_i=0$ for each $i=1,2,\ldots, n$. If $S$ is a graded *-ring, then  $S$ is \textit{graded positive definite} if $\sum_{i=1}^n s_i s_i^* =0$ implies that $s_i=0$ for homogeneous elements $s_i$.
\end{definition}



\begin{lem} \label{lem: pos-def}
Let $\G$ be an ultragraph. Then, $L_R(\G)$ is positive definite if and only if $R$ is positive definite.
\end{lem}
\begin{proof}
    Since $R$ is commutative and semisimple, $R=\Pi_i^n F_i$, where each $F_i$ is a field. Then $R$ is positive definite if and only if each $F_i$ is positive definite. By \cite[Proposition 3.4]{PRV12}, $L_{F_i}(E)$ is positive definite if and only if $F_i$ is positive definite. Hence, Proposition~\ref{prop:direct_prod_iso} implies that $L_R(E) = \Pi_i^n L_{F_i}(E)$ is positive definite if and only if each $F_i$ is positive definite.  
    
    Assume that $R$ is positive definite. Start by assuming that $\G$ has no singular vertices and let $E$ be its associated graph. By \cite[Theorem 9.5]{Fir20} there is a projection $Q$ in the multiplier algebra of $L_F(E)$ such that $L_F(\G)$ is isomorphic to the corner $QL_F(E)Q$. Thus, since $R$ is positive definite, $L_F(\G)$ is isomorphic to a subalgebra of the positive definite algebra $L_F(E)$, and, therefore, $L_R(\G)$ is positive definite.  

   Suppose now that $\G$ is an arbitrary ultragraph (possibly with singular vertices) and $R$ is positive definite. Then, the desingularization $\F$ of $\G$ is an ultragraph with no singular vertices, such that there is an injective *-homomorphism $\phi$ of $L_R(\G)$ into $L_R(\F)$ (\cite[Proposition 10.4.1]{Fir20}). Let $\{x_1, \ldots,x_m\} \subset L_R(\G)$, and suppose that $ x_1x_1^* + \cdots + x_mx_m^* =0$. Then,
    \begin{align*}
        0 &=  \phi(x_1x_1^*) + \cdots + \phi(x_mx_m^*) \\
        & =  \phi(x_1)\phi(x_1)^* + \cdots + \phi(x_m)\phi(x_m)^*. \\
    \end{align*}
    By the first part of the proof, $L_R(\F)$ is positive definite since $R$ is. Thus, $\phi(x_i) = 0$ for each $1\leq i\leq m$. Then, the injectivity of $\phi$ implies that $x_i=0$ for each  $1\leq i\leq m$, and, therefore, $L_R(\G)$ is positive definite. 

     Conversely, assume $L_R(\G)$ is positive definite and $r_1,\ldots,r_n\in R$ is such that $\sum_{i=1}^n r_ir_i* =0$. Choose any vertex $v\in G^0$. Then, since $p_v\in L_R(\G)$ is a projection,
    \begin{equation}\label{eq:pd}
    0 = \sum_{i=1}^n (r_ir_i^*)p_v = \sum_{i=1}^n (r_ip_v)(r_ip_v)^*.  
    \end{equation}
    Since $L_R(\G)$ is positive definite, it follows that $r_ip_v = 0$. Since $\{v\}\in \G^0$ is non-empty, it follows form \cite[Theorem 2.10]{ImaPouLar20} that $r_i=0$ for $i=1,2,\ldots, n$. Hence $R$ is positive definite, and the proof is complete. 
\end{proof}

We now prove our first main theorem, which characterizes ultragraph Leavitt path algebras with the various Rickart properties introduced in Section~\ref{s:prelim}.

\begin{thm} \label{thm: Rickart}
Let $\mathcal{G}$ be an ultragraph and $R$ a commutative semi-simple ring with unity. Then $L_R(\GG)$ is locally Rickart and graded locally Rickart. If, in addition, $R$ is positive definite, then $L_R(\GG)$ is a graded locally Rickart $*$-ring. Moreover, the following are equivalent.
\begin{enumerate}[(i)]
    \item $G^0\in \mathcal{G}^0$
    \item $L_R(\mathcal{G})$ is a Rickart ring. 
    \item $L_R(\mathcal{G})$ is a graded Rickart ring. 
\end{enumerate}
If $R$ is positive definite,  then (i)-(iii) are equivalent to
\begin{enumerate}[(iv)]
    \item $L_R(\mathcal{G})$ is a graded Rickart *-ring.
\end{enumerate}

\end{thm}
\begin{proof}
Let $R=\Pi_i^n F_i$, where each $F_i$ is a field. Let $E$ be the graph associated with $\mathcal{G}$ for an arbitrary enumeration of $\G^1$, as constructed in Section~\ref{s:prelim}.

By Proposition~\ref{prop:hereditary}, $L_R(\G)$ is hereditary and thus semihereditary. By \cite[Lemma 8]{HaVa18},  $L_{R}(\G)$ is locally Rickart. Every graded ring which is locally Rickart is graded locally Rickart. Hence, $L_R(\G)$ is graded locally Rickart.

If $R$ is positive definite, then $L_R(\G)$ is positive definite by Lemma~\ref{lem: pos-def}, which in turn implies that $L_{F_i}(\G)$ is positive definite. Hence, by \cite[Theorem 4.9]{HazNam23}, each $L_{F_i}(\G)$ is a graded von Neumann regular ring. By Corollary~\ref{cor:graded_regular}, $L_R(G)$ is graded von Neumann regular, and thus a graded locally Rickart $*$-ring, by \cite[Proposition 10]{HaVa18}. 

(i)$\Rightarrow$(ii): Suppose that $G^0\in\mathcal{G}^0$. Then $L_R(\GG)$ is a unital locally (graded) Rickart ring and, therefore, a Rickart ring. 

(ii) $\Rightarrow$ (iii): Follows from \cite[Lemma 2]{HaVa18}.

(iii) $\Rightarrow$ (i): Since $L_R(\G)$ is unital if and only if $G^0\in \G^0$, and (iii) implies 
 $L_R(\G)$ is unital, it follows that (iii) implies (i). Hence (i), (ii), and (iii) are equivalent.

Finally, if $R$ is positive definite, then $L_R(\G)$ is a locally graded Rickart $*$-ring by the first part of the proof. And, locally graded Rickart $*$-ring is a graded Rickart $*$-ring if and only if it is unital. Hence, (i) and (iv) are equivalent.
\end{proof}

\section{Baer ultragraph Leavitt path algebras} \label{s:baer}
In this section, we prove that every ultragraph Leavitt path algebra is locally Baer and we characterize when an ultragraph Leavitt path algebra is a Baer, graded Baer, and Baer *-ring. We need the following lemma for our main result.

\begin{lem} \label{lem: indistinguishability}
    Let $\G$ ultragraph, let $\{X_1, \ldots,X_n\}\subset \G^0$ be any finite collection of generalized vertices, and let $S\subset G^0 $. Let $A \in \G^0$ such that $A\cap S$ and $A \setminus S$ are infinite subsets of $G^0$. Suppose $A \cap S \notin \G^0$. Then, there exist two vertices $s \in A \cap S$ and $t\in A \setminus S$ such that $s\in X_i$ if and only if $t \in X_i$.
\end{lem} 
\begin{proof}
    For each element $v\in A$, associate a word $\omega(v) \in \{0,1\}^n$ having $i$th coordinate $\omega_i = 1$ if $v\in X_i$ and $0$ otherwise. Then there are at most $2^n$ words that partition the set $A$ into disjoint subsets. Define $r(\omega)=\{ v \in A : \omega= \omega(v)\}.$  First, notice that
    \begin{equation*}
        r(\omega)= \bigcap \limits_{i=1}^n \begin{cases} A \cap X_i \text{ if } \omega_i=1 \\
   A \setminus X_i \text{ if } \omega_i=0
    \end{cases}
    \end{equation*}    
    is a generalized vertex.
    
    We claim there must be a word $\omega$ such that $r(\omega)$ intersects both $S$ and its complement. Suppose not. Then, $r(\omega) \subset A \cap S$ or $r(\omega) \subset A \setminus S$, for each $\omega$. Therefore, 
    \[\bigcup_{r(\omega)\subset S}r(\omega) = A\cap S,\]
    which implies that $A\cap S$ is a generalized vertex (since there are at most $2^n$ words). This is a contradiction. Thus, there exist two elements $s\in A \cap S$ and $t \in A \setminus S$ such that $\omega(s)= \omega(t)$, so  $s\in X_i$ if and only if $t \in X_i$.
\end{proof}

\begin{thm} \label{thm: Baer}
Let $\mathcal{G}$ be an ultragraph and $R$ a commutative semi-simple ring with unity. Then the following are equivalent.
\begin{enumerate}[(i)]
    \item $\G$ is a row-finite, no-exit ultragraph in which every infinite path ends in a sink or a cycle
    \item $L_R(\mathcal{G})$ is a locally Baer ring. 
    \item $L_R(\mathcal{G})$ is a graded locally Baer ring. 
\end{enumerate}

If $R$ is positive definite, then (i)-(iii) are equivalent to

\begin{enumerate}[(iv)]
    \item $L_R(\mathcal{G})$ is a graded locally Baer *-ring. 
\end{enumerate}

As a corollary, the following conditions are equivalent:
\begin{enumerate}[(i')]
    \item $\G$ is a finite no-exit ultragraph where every infinite path ends in a sink or a cycle.
    \item $L_R(\mathcal{G})$ is a Baer ring. 
    \item $L_R(\mathcal{G})$ is a graded Baer ring. 
\end{enumerate}
If $R$ is positive definite, then (i')-(iii') are equivalent to 
\begin{enumerate}[(iv')]
    \item $L_R(\mathcal{G})$ is a graded Baer *-ring. 
\end{enumerate}

\end{thm}
\begin{proof}
(i) $\Rightarrow$ (ii) and (i) $\Rightarrow$ (iv): Since $R$ is semi-simple, it is isomorphic to a finite direct product  $\Pi_i^n F_i$, where each $F_i$ is a field. Let $E$ be the graph associated with $\mathcal{G}$, as constructed in Section~\ref{construction}. Since $\G$ is row-finite, $\Delta=\emptyset$, by Proposition~\ref{prop:Str_row_finite_equiv}. Then it follows from \cite[Proposition 3.23.1]{Fir20} that $L_R(\mathcal{G})\cong L_R(E)$. By Lemma~\ref{lem: inf path sink/cycles}, $E$ is a row-finite, no-exit graph in which every infinite path ends in a sink or cycle. Thus, by \cite[Theorem 15]{HaVa18}, each $L_{F_i}(E)$ is locally Baer, so as is $L_R(E)$ from Lemma~\ref{l:Rick_Baer_direct_prod}.  Therefore, so is $L_R(\G)$. If each $F_i$ is positive definite, then $R$ is positive definite by Lemma~\ref{lem: pos-def}, and thus, each $L_{F_i}(E)$ is a locally graded Baer *-ring by \cite[Theorem 15]{HaVa18}. Then $L_R(E)$ is a locally graded Baer *-ring by Lemma~\ref{l:Rick_Baer_direct_prod}. Therefore, so is $L_R(\G)$.

(ii) $\Rightarrow$ (iii): Any graded ring that is locally Baer is graded locally Baer.

(iii) $\Rightarrow$ (i): Suppose that $L_R(\G)$ is graded locally Baer. We claim that $\Delta=\emptyset$. Suppose not. Then, there is $n\in\N$ and $\omega\in \Delta_n$, which implies there is $i\leq n$ such that $|r(e_i)|=\infty$, where $e_i\in \G^1$. Put $A=r(e_i)$. Then, there exists an infinite subset $S\subset A$ such that $S$ is not a generalized vertex. This is because the set of infinite subsets of $A$ is uncountable, whereas the set of generalized vertices is countable. Notice that although $S$ is not a generalized vertex, we have that $p_s\in p_A L_R(\G) p_A$ for each $s\in S$. Hence, we may consider the left annihilator of the set of homogenous elements $P_S=\{p_s: s\in S\}$ in $p_A L_R(\G) p_A$. We claim that the left annihilator of the set $P_S$ is not generated as a left ideal by any idempotent.

Let $e=p_A x p_A$ be the homogenous idempotent generating $\ann_l(P_S)$. Then $x = \sum_{i=1}^n r_i s_{\alpha_i} p_{A_i} s_{\beta_i}^*$, where $r_i\in R$, $\alpha_{i}, \beta_i\in \G^*$, $A_i\in \G^0$  and $r(\alpha_i)\cap A_i\cap r(\beta_i)\neq \emptyset$ for $1\leq i\leq n$, by \cite[Theorem 2.9]{ImaPouLar20}. 
For $\mu=e_1e_2\ldots e_k$ with $e_i\in\G^1$, let $\mathrm{ver}(\mu)=\{v\in G^0\mid v=s(e_i)\,\,\text{ for }\,\, i=1\ldots,k\}$. Let 
$$\tilde{A} = A\setminus \{ v \in G^0: |\alpha_i|>0, \,\, |\beta_i|>0, \text{ and } v= s(\alpha_i) \text{ or } v\in \mathrm{ver}(\beta_i)\}.$$ 
Since $S \notin \G^0$, and we can recover $S$ from $\tilde{A} \cap S$ by adding a finite number of elements, $\tilde{A} \cap S \notin \G^0$. Therefore, both $ \tilde{A} \cap S$ and $\tilde{A} \setminus S$ are infinite sets.  Then, by Lemma~\ref{lem: indistinguishability}, we can pick two vertices $s \in \tilde{A} \cap S$ and $t\in \tilde{A} \setminus S$ such that for any set $X$ in the Boolean algebra generated by $\{ r(\alpha_i), A_i, r(\beta_i): i=1,\ldots,n\}$, we have that  $s\in X$ if and only if $t\in X$. Since $t \in \tilde{A} \setminus S$, it follows that $p_t\in \ann_l(P_S)$. Thus, $p_t = y e$ for some $y \in p_A L_R(G) p_A$. Then, $p_t = ye = ye^2= p_t e $, that is,  

\begin{eqnarray}\label{eq:6.3}
    p_t= p_te &=& p_tp_A\left( \sum_{i:|\alpha_i|=0}^nr_i s_{\alpha_i} p_{A_i} s_{\beta_i}^*  + \sum_{i:|\alpha_i|>0}^n r_i s_{\alpha_i} p_{A_i} s_{\beta_i}^*\right)p_A \nonumber  \\
    &=& \sum_{\substack{|\alpha_i|=0, \\ t\in r(\alpha_i)\cap A_i\cap r(\beta_i)}}^n r_i p_t s_{\beta_i}^* \nonumber \\
    &=& \sum_{\substack{|\alpha_i|=0,\,\, |\beta_i|=0 \\  t\in r(\alpha_i)\cap A_i\cap r(\beta_i)}}^n r_i p_t s_{\beta_i}^* + \sum_{\substack{|\alpha_i|=0,\,\, |\beta|>0 \\  t\in r(\alpha_i)\cap A_i\cap r(\beta_i)}}^n r_i p_t s_{\beta_i}^*  \nonumber \\
    &=& \sum_{\substack{|\alpha_i|=0,\,\, |\beta_i|=0 \\  t\in r(\alpha_i)\cap A_i\cap r(\beta_i)}}^n r_i p_t + \sum_{\substack{|\alpha_i|=0,\,\, |\beta|>0 \\  t\in r(\alpha_i)\cap A_i\cap r(\beta_i)}}^n r_i p_t s_{\beta_i}^*
\end{eqnarray}
To ease the notational burden, put 
$$u = \displaystyle \sum_{\substack{|\alpha_i|=0,\,\, |\beta_i|=0 \\  t\in r(\alpha_i)\cap A_i\cap r(\beta_i)}}^n r_i p_t \,\,\, \text{ and }\,\,\,  v = \sum_{\substack{|\alpha_i|=0,\,\, |\beta|>0 \\  t\in r(\alpha_i)\cap A_i\cap r(\beta_i)}}^n r_i p_t s_{\beta_i}^*, $$
so that $p_t-u = v$. Since $\mathrm{ver}(\beta_i)\notin \tilde{A}$ for each $1\leq i\leq n$, it follows that
\[(p_t-u)p_t = vp_t =0.\]
Hence $(p_t-u)p_t = p_t - u = 0$, which implies that $v=0$. Reorder and relabel the terms of $v$ so that 
\begin{equation}\label{eq:og_v}
    v = \sum_{i=1}^m r_ip_ts_{\beta_i}^*,
\end{equation}
with $|\beta_1|\geq |\beta_2|\cdots\geq |\beta_m|$ and $\beta_i\neq \beta_j$ for each $1\leq i,j\leq m$. We claim that $r_i=0$ for $i=1,\ldots,m$. We show this iteratively. For step one, multiply $v$ by $s_{\beta_1}$ from the right. Then, only summands $r_jp_ts_{\beta_j}^*s_{\beta_1}$ for which $|\beta_1|>|\beta_j|$ and $\beta_1 = \beta_j\gamma_j$ for some $\gamma_j\in \G^*$ with $|\gamma_j|\geq 1$, are nonzero (see for example \cite{Tom03}). Hence, multiplying $v$ by $s_{\beta_1}$ from the right yields a sum of the form
\[0 = r_1p_tp_{r(\beta_1)} +\sum_{j=2}^{m_1}r_jp_ts_{\gamma_j},\]
with $m_1\leq m$. Since $s(\gamma_j)\notin \tilde{A}$, the sum $\sum_{j=2}^{m_1}r_jp_ts_{\gamma_j}=0$. Thus, $ 0=r_1p_tp_{r(\beta_1)} = r_1p_t$. Since $p_t\neq 0$, \cite[Theorem 2.10]{ImaPouLar20} implies that $r_1=0$. 

For step two, multiply $v$ by $s_{\beta_2}$ from the right. Then proceeding similarly to the first step, we see that $r_2 = 0$. Continuing in this way, after a finite number of iterations, it follows that $r_1=\cdots=r_m =0$. Hence, 
\begin{equation}\label{eq:6.5}
    \sum_{\substack{|\alpha_i|=0,\,\, |\beta|>0 \\  t\in r(\alpha_i)\cap A_i\cap r(\beta_i)}}^n r_i s_{\beta_i}^* = 0,
\end{equation}
and then Equation~(\ref{eq:6.3}) implies that 
\begin{eqnarray} \label{eq:6.6}
    \sum_{\substack{|\alpha_i|=0,\,\, |\beta_i|=0 \\  t\in r(\alpha_i)\cap A_i\cap r(\beta_i)}}^n r_i  = 1_R.
\end{eqnarray}



Now, 
\begin{eqnarray}
      p_se  &=& p_sp_A\left( \sum_{i:|\alpha_i|=0}^nr_i s_{\alpha_i} p_{A_i} s_{\beta_i}^*  + \sum_{i:|\alpha_i|>0}^n r_i s_{\alpha_i} p_{A_i} s_{\beta_i}^*\right)p_A \nonumber \\
      &=& \sum_{\substack{|\alpha_i|=0,\,\, |\beta_i|=0 \\  s\in r(\alpha_i)\cap A_i\cap r(\beta_i)}}^n r_i p_s + \sum_{\substack{|\alpha_i|=0,\,\, |\beta|>0 \\  s\in r(\alpha_i)\cap A_i\cap r(\beta_i)}}^n r_i p_s s_{\beta_i}^*,\label{eq:20}
\end{eqnarray}
where the last equality follows from the fact that $s\in \tilde{A}$.
Since $\omega(s)=\omega(t)$ (as denoted in the proof of Lemma~\ref{lem: indistinguishability}), $s\in r(\alpha_i)\cap A_i \cap r(\beta_i)$ if and only if $t\in r(\alpha_i)\cap A_i \cap r(\beta_i)$. Thus, Equation~(\ref{eq:20}) may be rewritten as
\begin{eqnarray*}
     p_se  &=& p_s \left(\sum_{\substack{|\alpha_i|=0,\,\, |\beta_i|=0 \\  s\in r(\alpha_i)\cap A_i\cap r(\beta_i)}}^n r_i\right)  + p_s \left( \sum_{\substack{|\alpha_i|=0,\,\, |\beta|>0 \\  s\in r(\alpha_i)\cap A_i\cap r(\beta_i)}}^n r_i s_{\beta_i}^*\right) \\
     &=& p_s \left(\sum_{\substack{|\alpha_i|=0,\,\, |\beta_i|=0 \\  t\in r(\alpha_i)\cap A_i\cap r(\beta_i)}}^n r_i\right)  + p_s \left( \sum_{\substack{|\alpha_i|=0,\,\, |\beta|>0 \\  t\in r(\alpha_i)\cap A_i\cap r(\beta_i)}}^n r_i s_{\beta_i}^*\right) \\
    &=&  p_s,
\end{eqnarray*}
where the last equality follows from Equations~(\ref{eq:6.5}) and (\ref{eq:6.6}).

Since $p_s=p_Ap_sp_A \in p_A L_R(\G) p_A$ and $p_se=p_s$, it follows that 
\[p_s\in (p_AL_R(\G)p_A)e= \ann_l^{p_AL_R(\G)p_A}(P_S),\]
which is a contradiction since $s\in S$. Hence, $\Delta=\emptyset$, and 
therefore, $L_R(\G)\cong L_R(E)$, by \cite[Proposition 3.23.1]{Fir20}. Hence, $L_R(E)$ is locally Baer. Since $L_R(E)$ is locally Baer, each $L_{F_i}(E)$ is locally Baer by Lemma \ref{l:Rick_Baer_direct_prod}. Therefore, by \cite[Theorem 15]{HaVa18}, $E$ must be a row-finite no-exit graph in which every infinite path ends in a sink or cycle. Therefore, $\G$ is a row-finite, no-exit ultragraph in which every infinite path ends in a sink or a cycle, by Lemma~\ref{lem: inf path sink/cycles}.

(iv) $\Rightarrow$ (iii): Any locally graded Baer *-ring is a locally graded Baer ring since every projection is an idempotent.

The equivalence of (i')-(iv') follows form the equivalence of (i)-(iv) when $L_R(\G)$ is unital. 
\end{proof}

\begin{remark}
    In Theorem~\ref{thm:baer*} below we extend \cite[Theorem 16]{HaVa18} to ultragraphs. Note that there is a typographical error in \cite[Theorem 16]{HaVa18}. In addition to the assumption that the graph $E$ must be a disjoint union of graphs that are acyclic or isolated loops, all acyclic graphs must also be row-finite, as is evident from the assumptions in \cite[Theorem 15]{HaVa18} and the proof of \cite[Theorem 16]{HaVa18}. 
\end{remark}

\begin{thm}\label{thm:baer*}
Let $\mathcal{G}$ be an ultragraph and $R$ a positive definite commutative semi-simple ring with unity. Then the following are equivalent.
\begin{enumerate}[(i)]
    \item $\G$ is a disjoint union of ultragraphs that are either acyclic and row-finite with each infinite path ending in a sink, or isolated loops.
    \item $L_R(\mathcal{G})$ is a locally Baer *-ring. 
\end{enumerate}

As a corollary, the following conditions are equivalent:
\begin{enumerate}[(i')]
    \item $\G$ is a finite disjoint union of ultragraphs that are finite and acyclic, or isolated loops.
    \item $L_R(\mathcal{G})$ is a Baer *-ring. 
\end{enumerate}
\end{thm}

\begin{proof}
(i) $\Rightarrow$ (ii): Suppose $\G$ is a disjoint union of ultragraphs that are row-finite and acyclic, or isolated loops.  Let $E$ be the graph associated with $\mathcal{G}$, as constructed in Section~\ref{construction}. Then $\Delta=\emptyset$, by Proposition~\ref{prop:Str_row_finite_equiv}. Since $R$ is semi-simple, it is isomorphic to a finite direct product  $\Pi_i^n F_i$, where each $F_i$ is a field, and by Proposition \ref{prop:direct_prod_iso}, $L_R(\G) \cong \Pi_i^n L_{F_i} (\G)$.  Since $\Delta$ is empty, \cite[Proposition 3.23.1]{Fir20} implies that $L_{F_i}(\mathcal{G})\cong L_{F_i}(E)$. Then $E$ is a disjoint union of graphs that are acyclic and row-finite, or isolated loops, by Theorem \ref{lem: inf path sink/cycles}. Therefore, each $L_{F_i}(E)$ is a locally Baer *-ring, by \cite[Proposition 16]{HaVa18}, and therefore, so are $L_{F_i}(\G)$ and $L_R(\G)$, by Lemma \ref{l:Rick_Baer_direct_prod}.

(ii) $\Rightarrow$ (i): Suppose $L_R(\G)$ is a locally Baer *-ring. Since a locally Baer *-ring is locally Baer, Theorem \ref{thm: Baer} and Proposition~\ref{prop:Str_row_finite_equiv} imply that $\Delta$ is empty. Thus, similarly to the first part of the proof, this implies that each $L_{F_i}(E)$ is a locally Baer *-ring. From \cite[Proposition 16]{HaVa18}, this occurs precisely when $E$ is a disjoint union of graphs that are row-finite and acyclic, or isolated loops. Since $\Delta$ is empty, all paths in $\G$ are in bijective correspondence with paths in $E$ of the same length. Thus, $\G$ is a finite disjoint union of ultragraphs that are row-finite and acyclic, or isolated loops.

(i') $\Rightarrow$ (ii'): If this is the finite disjoint union of finite ultragraphs, then (i) holds and $L_R(\G)$ is unital. Therefore, by looking at the corner generated by the unit, $L_R(\G)$ is a Baer *-ring.

(ii') $\Rightarrow$ (i'): Suppose $L_R(\G)$ is a Baer *-ring. Thus, $L_R(\G)$ is unital and a locally Baer *-ring, which implies $\G$ is a disjoint union of ultragraphs that are row-finite and acyclic or isolated loops, since (ii) implies (i). Since each of the ultragraphs are row-finite, only finite collections of vertices are elements in $\G^0$. Therefore, for $L_R(\G)$ to be unital, $\G$ can only contain finitely many vertices, so both the disjoint union is finite, and each ultragraph is finite as well. If either were infinite, we would need $\G^0$ to be closed under infinite unions, which is not required.
\end{proof}

\bibliographystyle{abbrv}
\bibliography{ref}
\end{document}